\documentclass[12pt,a4paper,fleqn]{article}
\usepackage{exscale,xypic}
\usepackage[intlimits]{amsmath} \topmargin-2cm \textheight242mm
\newfont{\tabfont}{cmr8 at 10pt}

\begin{document}
\begin{center}
{\Large \bf Open-multicommutativity of some functors related to the functor of
probability measures}
\end{center}

\begin{center}
{\bf Roman Kozhan}\\Department of Mechanics and Mathematics\\ Lviv
National University
\end{center}


\begin{center}
Abstract
\end{center}

\vspace{0.2cm}{\tabfont The property of a normal functor to be
open-multicommutative proposed by Kozhan and Zarichnyi (2004) is investigated. A
number of normal functors related to the functor of probability measures and
equipped with convex structure are considered here and it is proved that the
functors $cc$, $ccP$, $G_{cc}P$ and $\lambda_{cc}P$ are open-multicommutative.}

\begin{center}
{\large \bf 1. Introduction}
\end{center}

The classical object in topology and functional analysis - spaces of probability
measures - is widely used in economics and game theory last ten years (see Lucas
(1971), and Prescott (1971), Mas-Colell (1984), Jovanovi\'c and Rosenthal
(1988)). All this investigations deal with the notion of the set-valued
correspondence map which assigns to every probability measures on the factors of
the product of compacta the set of probability measures with these marginals,
\[\psi\left(\mu_1,...,\mu_n\right)=\{\lambda\in
P\left(X_1\times...\times X_n\right)\colon
P\pi_i\left(\lambda\right)=\mu_i,\text{  }i=1,...,n\}.
\]

The problem of continuity of this map can be equivalently
redefined in terms of openness of the characteristic map of the
 bicommutative diagram
\[\xymatrix{P\left(X\times Y\right)\ar[rr]^{P\pi_1}\ar[d]^{P\pi_2}&&P\left(X\right)
\ar[d]^{P{\mathbf 1}_{*}}\\P\left(Y\right) \ar[rr]_{P{\mathbf
1}_{*}}&&P\left(\{*\}\right)}
\]

Due to the well known theorem of Ditor and Eifler (1972) the
openness property of the characteristic map is closely related to
the property of the bicommutativity of a normal functor and, in
particular, the functor of probability measures.

A common generalization of these two properties generates a new notion of
multi-commutativity of the normal function which was proposed by Kozhan and
Zarichnyi (2004). In this paper they investigate the functor of probability
measures and show that this functor is multi-commutative.

In the economic theory, together with the space of probability measures,
different convex structures are commonly used, in particular, the spaces of
convex closed subsets of the space of probability measures have very nice
application.

Here the open-multicommutativity property of normal functors is investigated for
functors $cc$, $ccP$, $G_{cc}P$ and $\lambda_{cc}P$.

\begin{center}
{\large \bf 2. Notations and definitions}
\end{center}

Let ${\mathbf K}$ be a finite category. Denote by $|{\mathbf K}|$
the class of all objects of the category ${\mathbf K}$. For every
$A,B\in |{\mathbf K}|$ the set ${\mathbf K}\left(A,B\right)$
consists of all morphisms from $A$ to $B$ in ${\mathbf K}$. A
functor $D\colon{\mathbf K}\rightarrow {\mathbf {Comp}}$ is called
a {\em diagram}.

{\bf Definition 1.} The set of morphisms
\begin{equation}
\left( X \overset{g_A}{\rightarrow} D\left( A\right)\right)_{A \in
|{\mathbf K}| }
\end{equation}
is said to be a {\em cone} over the diagram $D$ if and only if for
every objects $A,B \in |{\bf K}|$ and for every morphism
$\varphi\colon A \rightarrow B$ in ${\mathbf K}$ the diagram
\[\xymatrix{&X\ar[dl]_{g_A}\ar[dr]^{g_B}&\\D\left(A\right)
\ar[rr]_{D\left(\varphi\right)}&&D\left(B\right)}
\]is commutative.

{\bf Definition 2.} Cone (1) is called a {\em limit} of the diagram $D$ if the
following condition is satisfied: for each cone $C^\prime=\left( X^\prime
\overset{{g^\prime_A}}{\rightarrow} D\left( A\right)\right)_{A \in |{\mathbf K}|
}$ there exists a unique morphism $\chi_{C^\prime}\colon X^\prime \rightarrow X$
such that $g^\prime_A = g_A \circ \chi_{C^\prime}$ for every $A
\in |{\mathbf K }|$.

Further, we denote this cone by $\lim\left(D\right)$. The map $\chi_{C^\prime}$
is called the {\em characteristic map} of $C^\prime$.

{\bf Definition 3.} The cone $C^\prime=\left( X^\prime
\overset{{g^\prime_A}}{\rightarrow} D\left( A\right)\right)_{A \in
|{\mathbf K}| }$ is called {\em open-multicommutative} if the
characteristic map $\chi_{C^\prime}$ is open and surjective.

Let {\bf F} be a normal functor in the category {\bf {Comp}}. Define the diagram
${\mathbf F}\left(D\right)\colon{\mathbf K}\rightarrow {\mathbf {Comp}}$ in the
following way: for every $A \in {\mathbf K}$ let ${\mathbf
F}\left(D\right)\left(A\right)={\mathbf F}\left(D\left(A\right)\right) $ and for
every morphism $\varphi\in{\mathbf K}\left(A,B\right)_{A,B\in |{\mathbf K}|}$ we
put $ {\mathbf F}\left(D\right)\left(\varphi\right)={\mathbf
F}\left(D\left(\varphi\right)\right)$.

{\bf Definition 4.} The normal functor ${\bf F}$ is called {\em
open-multicommutative} if it preserves open-multicommutative
cones, i.e. the cone
\[{\mathbf F}\left( C^\prime \right)=\left( {\mathbf F}\left(
X^\prime \right) \overset{{\mathbf F}g^\prime_A}{\rightarrow}
{\mathbf F}\left( D\left( A\right) \right)\right)_{A \in |{\mathbf
K}| }
\]over the diagram ${\mathbf F}\left(D\right)$ is
open-multicommutative.

For the normal functor ${\mathbf F}$ we assume that the cone
$\left(Y_{\mathbf F }\overset{\pi_A}{\rightarrow}{\mathbf F
}\left(D\left(A\right)\right)\right)_{A\in |{\mathbf K}|}$ is a
limit of the diagram ${\mathbf F}\left(D\right)$.

Let us assume that $X$ is a compact subspace of some locally
convex space $E$. We consider a functor $cc\colon {\mathbf
{Conv}}\rightarrow {\mathbf {Comp}}$. It is defined as
\[cc\left(X\right)=\{A\subset X\colon A \text{ is closed and convex}\}\subset \exp\left(X\right).
\]

Recall the constructions of well-known inclusion hyperspace and
superextension functors. We set
\[G\left(X\right)=\{{\mathbf A}\in \exp^2\left(X\right)\colon A\in {\mathbf A}\text{ and
}A\subset B\in \exp\left(X\right)\Rightarrow B\in{\mathbf A}\}.
\]\[\lambda\left(X\right)=\{{\mathbf A}\in G\left(X\right)\colon {\mathbf A}\text{
is a maximal linked system }\}.
\] For properties of the functors $G$ and $\lambda$ reader is
referred to Teleiko and Zarichnyi (1999).

In order to combine these constructions with convex structure we use the defined
above functors and the functor $cc$. In such a way we can define a space
\[G_{cc}\left(X\right)=\{{\mathbf A}\subset cc\left(X\right)\colon\text{
}{\mathbf A}\text{ is closed and }A\in {\mathbf A}\text{ and
}A\subset B\in cc\left(X\right)\Rightarrow B\in{\mathbf A}\}.
\]It is obviously that $G_{cc}\left(X\right)\subset G\left(X\right)$
and we can endow it with the topology induced from $G\left(X\right)$. For a map
$f\colon X\rightarrow Y$ in ${\mathbf {Comp}}$ and ${\mathbf A}\in
G_{cc}\left(X\right)$ let
\[G_{cc}f\left({\mathbf A}\right)=\{B\in G_{cc}\left(Y\right)|\text{
}f\left(A\right)\subset B\text{, }A\in {\mathbf A}\}.
\] Thus we can define a covariant functor $G_{cc}\colon {\mathbf {Conv}} \rightarrow
{\mathbf {Comp}}$. In the same way we can define a functor
$\lambda_{cc}\colon {\mathbf {Conv}}\rightarrow {\mathbf {Comp}}$.

Due to the natural convex structure of the functor $P$ we can
compose the functors $G_{cc}$ and $\lambda_{cc}$ with the functor
of probability measures and it gives us two new constructions of
functors in the category ${\mathbf {Comp}}$:
\[G_{cc}P=G_{cc}\circ P\colon {\mathbf {Comp}}\rightarrow {\mathbf {Comp}}
\]\[\lambda_{cc}P=\lambda_{cc}\circ P\colon {\mathbf
{Comp}}\rightarrow {\mathbf {Comp}}.
\]

\begin{center}
{\large \bf 3. Open-multicommutativity of the functor $ccP$ and
related functors}
\end{center}

{\bf Proposition 1.} {\it A normal bicommutative functor ${\mathbf
F}$ is multicommutative.}

{\bf Proof.} This Proposition is actually proved in Kozhan and Zarichnyi (2004)
for the functor of probability measures. However, in the proof they use only
normality and bicommutativity of the functor $P$ and that is why it can be
applied for every normal bicommutative functor.

Q.E.D.

\medskip

The following proposition is also the generalization of the result
of Kozhan and Zarichnyi.

{\bf Proposition 2.} {\it A normal open functor ${\mathbf F}$ is
open-multicommutative if and only if the characteristic map
\[\chi_{\mathbf F}\colon {\mathbf F}\left(X\right)\rightarrow Y_{\mathbf F}
\] is open for every diagram $D$ with finite spaces $D\left(A\right)$, $A\in |{\mathbf K}|$.}

{\bf Proof.} The proof is analogical to the proof of the theorem 1
in Kozhan and Zarichnyi (2004). The special properties of functor
$P$ are used only in the case of finite spaces $D\left(A\right)$.
The rest of the proof consider only general properties of normal
bicommutative functors and this implies that the scheme of the
proof is the same.

Q.E.D.

\medskip

{\bf Proposition 3.} {\it Functors $\exp$, $G$ and $\lambda$ are
open-multicommutative.}

{\bf Proof.} Let us check whether the conditions of the
proposition 2 are satisfied. The functor $\exp$ is normal and
bicommutative (see Teleiko and Zarichnyi (1999)). Assume that
every compactum $D\left(A\right)$ is finite for each $A\in
|{\mathbf K}|$. Since every finite compactum $D\left(A\right)$ is
discrete space therefore $\exp\left(D\left(A\right)\right)$ is
also discrete, which follows from the properties of the Vietoris'
topology. It is known (see Kozhan and Zarichnyi (2004)) that
$Y_{\exp}\subseteq \underset{A\in |{\mathbf
K}|}{\prod}\exp\left(D\left(A\right)\right)$, which is discrete as
well as $\exp\left(X\right)$, is a subset of the discrete space
$\exp\left(\underset{A\in |{\mathbf K}|}{\prod}D\left(A\right)
\right)$. Thus the characteristic map $\chi_{\exp}$ is the map of
two discrete spaces and this necessarily implies that it is open.
The spaces $G\left(X\right)$ and $\lambda\left(X\right)\subset
\exp^2\left(X\right)$ for every compactum $X$ and therefore are
also discrete if $X$ is so. Then normality and bicommutativity of
them (see Proposition 2) implies the open-multicommutativity of
these functors.

Q.E.D.

\medskip

{\bf Lemma 1.} {\it Let $B\subset T\times T$ and
$\pi_1\left(B\right)\subset C\subset T$. If $C$ is convex set then
this implies that
$\pi_1\left(\text{conv}\left(B\right)\right)\subset C$. If in
addition we have that $\pi_1\left(B\right)\supset C$ then
$\pi_1\left(\text{conv}\left(B\right)\right)=C$.}

{\bf Proof.} Consider an arbitrary point $x\in
\text{conv}\left(B\right)\setminus B$. There exist two points
$x_1, x_2\in B$ and $\alpha \in \left(0,1\right)$ such that $x=\alpha
x_1+\left(1-\alpha\right)x_2$. Since
$\pi_1\left(x_1\right),\pi_1\left(x_2\right)\in C$ and $C$ is convex,
$\pi_1\left(x\right)\in C$. The point $x$ is arbitrary chosen thus the first
statement of the lemma is proved.

The second statement of the lemma is evident.

Q.E.D.

\medskip

{\bf Proposition 4.} {\it The functor $cc \colon{\mathbf
{Conv}}\rightarrow {\mathbf {Comp}}$ is open-multicommutative.}

{\bf Proof.} Let us prove that the characteristic map $\chi_{cc}\colon
cc\left(X\right)\rightarrow Y_{cc}$ is open. Consider an arbitrary point $B\in
cc\left(X\right)$ and an arbitrary sequence $\{C_i\}_{i\in\aleph}\subset Y_{cc}$
such that $\lim C_i=C$, where $C=\chi_{cc}\left(B\right)$. Recall that for every
compactum $X$ with convex structure the space $cc\left(X\right)\subset
\exp\left(X\right)$. Since $B\in\exp\left(X\right)$ and $\chi_{\exp}$ is open map
this implies that there exists a sequence $\{B_i\}_{i\in\aleph}\subset
\exp\left(X\right)$ such that
\[\lim B_i=B \text{ and }\chi_{\exp}\left(B_i\right)=C_i
\]for every $i\in\aleph$. Denote
by $D_i$ the convex hull of the set $B_i$. Since the function $\text{conv}\colon
\exp\left(T\right)\rightarrow cc\left(T\right)$ is continuous in the Vietoris
topology for every compactum $T$ with the convex structure, we see that
\[\lim D_i=\lim \text{conv}\left(B_i\right)=\text{conv}\left(\lim
B_i\right)=\text{conv}\left(B\right)=B\in cc\left(X\right).
\] Each $C_i$ is a convex set therefore $\pi_A \left(C_i\right)$ is
also a convex set for every $A\in |{\mathbf K}|$ hence Lemma 1
implies that $\chi_{cc}\left(D_i\right)=C_i$ for every
$i\in\aleph$. Thus, the map $\chi_{cc}$ is open. The surjectivity
of the characteristic map is obvious since every element in
$\underset{A\in|{\mathbf K}|}{\prod} cc\left(D\left(A\right)
\right)$ is also in $cc\left(\underset{A\in|{\mathbf
K}|}{\prod}D\left(A\right) \right)$.

Q.E.D.

\medskip

{\bf Proposition 5.} {\it Let categories $Q_1$, $Q_2$, $Q_3\subset {\mathbf Top}$
and functors $F_1\colon Q_1\rightarrow Q_2$ and $F_2\colon Q_3\rightarrow Q_1$
are open-multicommutative then the composition $F_1\circ F_2\colon Q_3\rightarrow
Q_2$ is also open-multicommutative.}

{\bf Proof.} The open-multicommutativity of the functor $F_2$ implies that the
characteristic map $\chi_{F_2,D}\colon F_2\left(X\right) \rightarrow Y_{F_2}$ is
open and surjective. Since the functor $F_1$ is open hence the composition
$F_1\chi_{F_2,D}\colon F_1\circ F_2\left(X\right)\rightarrow F_1\left(
Y_{F_2}\right)$ is also open and surjective. Consider now a diagram $D_1$ in the
category ${\mathbf K}$ such that for every $A\in |{\mathbf K}|$ we have
$D_1\left(A\right)=F_2\left(D\left(A\right)\right)$ and for each $\varphi\in
{\mathbf K}\left(A,B\right)_{A,B\in|{\mathbf K}|}$ we have
$D_1\left(\varphi\right)=F_2\left(D_1
\left(\varphi\right)\right)$. The functor $F_1$ is
open-multicommutative, so this implies that the characteristic map
$\chi_{F_1,D_1}\colon F_1\left(Y_{F_2}\right)\rightarrow Y_{F_1\circ F_2}$ is
open and surjective. Thus, a map $\chi_{F_1\circ
F_2,D}=F_1\chi_{F_2,D}\circ\chi_{F_1,D_1}\colon F_1\circ F_2
\left(X\right)\rightarrow Y_{F_1\circ F_2}$ is open and surjective for any
diagram $D$ as the composition of two open and surjective maps. This implies that
the functor $F_1\circ F_2$ is open-multicommutative.

Q.E.D.

\medskip

{\bf Corollary 1.} {\it The functor $ccP\colon {\mathbf
{Comp}}\rightarrow{\mathbf {Comp}}$ is open-multicommutative.}

{\bf Proof.} This follows from the open-multicommutativity of the functors $P$
(see Kozhan and Zarichnyi (2004)), $cc$ (see Proposition 4) and Proposition 5.

Q.E.D.

\medskip

{\bf Proposition 6.} {\it The functors $G_{cc}$ and $\lambda_{cc}$
are open-multicommutative.}

{\bf Proof.} Let us define a retraction for every
\[r_{cc}X\colon G\left(X\right)\rightarrow G_{cc}\left(X\right)
\]in the following way: for every ${\mathbf A}\in G\left(X\right)$
\[r_{cc}X\left({\mathbf A}\right)=\{B\in
\exp\left(X\right)|B=\text{conv}\left(A\right)\text{, }A\in
{\mathbf A}\}.
\]It is easy to see that $r_{cc}X\left({\mathbf A}\right)\in
G_{cc}\left(X\right)$.

The base of the space $G\left(X\right)$ is formed by the sets $\left(U_1^+\cap
...\cap U_m^+\right)\cap\left(V_1^-\cap ...\cap V_n^+\right)$, where
\[U^+=\{{\mathbf A}\in G\left(X\right)|\text{ there exists
}A\in{\mathbf A}\text{ such that }A\subset U\}
\]
\[U^-=\{{\mathbf A}\in G\left(X\right)|\text{ }A\cap
U\neq\emptyset\text{ for every } A\in {\mathbf A}\}
\]for some open set $U\subset X$.

Let us show that $r_{cc}X$ is continuous for every $X\in {\mathbf {Conv}}$. To
prove the continuity at a point ${\mathbf A}_0\in G\left(X\right)$ it is
sufficient to show that for every element $U^+\cap G_{cc}\left(X\right)$
(respectively $U^-\cap G_{cc}\left(X\right)$) which contains ${\mathbf A}_0$
there exists a neighborhood ${\mathbf V}\subset G\left(X\right)$ of the point
${\mathbf B}_0=r_{cc}X\left({\mathbf A}_0\right)$ such that
\[r_{cc}X\left({\mathbf V}\right)\subset U^+\cap
G_{cc}\left(X\right)\text{ ( } U^-\cap G_{cc}\left(X\right)\text{
respectively) }.
\]Assume first that ${\mathbf U}=U^-\cap G_{cc}\left(X\right)$.
Denote ${\mathbf V}=U^-$. Then for every ${\mathbf A}\in{\mathbf
V}$ we have
\[\forall A\in{\mathbf A}\colon A\cap
U\neq\emptyset\Rightarrow\text{conv}\left(A\right)\cap U
\neq\emptyset.
\]Since
\[{\mathbf B}=r_{cc}X\left({\mathbf A}\right)=\{B\in
\exp\left(X\right|B=\text{conv}\left(A\right)\text{, }A\in
{\mathbf A}\},
\]this implies that $\forall B\in{\mathbf B}$ we have $B\cap U \neq\emptyset$
and then ${\mathbf B}\in U^-\cap G_{cc}\left(X\right)={\mathbf
U}$.

Assume now that ${\mathbf U}=U^+\cap G_{cc}\left(X\right)$. Since $X$ is a subset
of locally convex space $E$, there exists a base of $X$ consisting of convex
sets. This implies that we can find an open convex set $V\subset U$. Denote
${\mathbf V}=V^+$. For every ${\mathbf A}\in {\mathbf V}$ there is a set $A_1\in
{\mathbf A}$ such that $A_1\subset V$. The set $V$ is convex, thus
$\text{conv}\left(A_1\right)\subset V\subset U$. This implies that for ${\mathbf
B}=r_{cc}X\left({\mathbf A}\right)$ we can find $B_1\in {\mathbf B}$ such that
$B_1\subset U$ and this means that ${\mathbf B}\in U^+\cap
G_{cc}\left(X\right)={\mathbf U}$. Thus, the map $r_{cc}$ is continuous.

Let us prove that the characteristic map $\chi_{G_{cc}}\colon
G_{cc}\left(X\right)\rightarrow Y_{G_{cc}}$ is open. Consider an arbitrary point
${\mathbf B}\in G_{cc}\left(X\right)$ and an arbitrary net $\{{\mathbf
C_i}\}_{i\in\aleph}\subset Y_{G_{cc}}$ such that $\lim {\mathbf C_i}={\mathbf
C}$, where ${\mathbf C}=\chi_{G_{cc}}\left(B\right)$. For every compactum $X$ the
space $G_{cc}\left(X\right)\subset G\left(X\right)$ and then ${\mathbf B}\in
G\left(X\right)$. Since $\chi_{G}$ is open map this implies that there exists a
net $\{{\mathbf B_i}\}_{i\in\aleph}\subset G\left(X\right)$ such that
\[\lim {\mathbf B_i}={\mathbf B} \text{ and }\chi_{G}\left({\mathbf
B_i}\right)={\mathbf C_i}
\]for every $i\in\aleph$.

Denote by ${\mathbf D_i}$ an image of the function $r_{cc}$ of the set ${\mathbf
B_i}$. Since the function $r_{cc}$ is continuous, we see that
\[\lim {\mathbf D_i}=\lim r_{cc}\left({\mathbf
B_i}\right)=r_{cc}\left(\lim {\mathbf
B_i}\right)=r_{cc}\left({\mathbf B}\right)={\mathbf B}\in
G_{cc}\left(X\right).
\]Each ${\mathbf C_i}$ is in $\underset{A\in |{\mathbf
K}|}{\prod}G_{cc}\left(D\left(A\right)\right)$ and let ${\mathbf
C^{\prime}_i}=\chi_{G_{cc}}\left({\mathbf D_i}\right)$ for every $i\in \aleph$.
Let us prove that ${\mathbf C_i}={\mathbf C^{\prime}_i}$. For every $A\in
|{\mathbf K}|$ we see that
\[G_{cc}\pi_A\left({\mathbf D_i}\right)=\pi_A\left({\mathbf
C^{\prime}_i}\right)=\{C^{\prime}\in cc\left(X\right)|\text{
}\pi_A\left(D\right)\subset C^{\prime}\text{, }D\in {\mathbf
D_i}\}.
\]Lemma 1 proves that if for some $C\in cc\left(X\right)$ the
statement $\pi_A\left(B\right)\subset C$ is satisfied for $B\in {\mathbf B_i}$
then $\pi_A\left(D\right)=\pi_A\left(\text{conv}\left(B\right)\right)\subset C$.
This implies that $G_{cc}\pi_A\left({\mathbf B_i}\right)\subset
G_{cc}\pi_A\left({\mathbf D_i}\right)$. On the other hand, if
$\pi_A\left(D\right)\subset C^{\prime}$ for some $D\in {\mathbf D_i}$ then
\[\pi_A\left(B\right)\subset
\pi_A\left(\text{conv}\left(B\right)\right)=\pi_A\left(D\right)\subset
C^{\prime}$ for $B\in{\mathbf B_i}.
\]This means that
$G_{cc}\pi_A\left({\mathbf B_i}\right)\supset G_{cc}\pi_A
\left({\mathbf D_i}\right)$. Since these two inclusions are
satisfied for every $A\in |{\mathbf K}|$, we have
\[\underset{A\in |{\mathbf K}|}{\prod}G_{cc}\pi_A\left({\mathbf
B_i}\right)=\underset{A\in |{\mathbf
K}|}{\prod}G_{cc}\pi_A\left({\mathbf D_i}\right)
\]and it is equivalent to ${\mathbf C_i}={\mathbf C^{\prime}_i}$. This
immediately implies that $\chi_{G_{cc}}\left({\mathbf D_i}\right)={\mathbf C_i}$
for every $i\in \aleph$. Thus, the characteristic map $\chi_{G_{cc}}$ is open and
the functor $G_{cc}$ is open-multicommutative. This map is also surjective since
for every ${\mathbf C}\in Y_{G_{cc}}$ we have $\chi_{G_{cc}}\left({\mathbf
D_C}\right)={\mathbf C}$, where
\[{\mathbf D_C}=\{D\in cc\left(X\right)|\text{ }D\supset
\underset{A\in |{\mathbf K}|}{\prod}C_A\text{, }C_A\in
\pi_A\left({\mathbf C}\right)\}.
\]

Let us note that the restriction
\[r_{cc}|_{\lambda\left(X\right)}\colon
\lambda\left(X\right)\rightarrow \lambda_{cc}\left(X\right)
\]is also continuous retraction. Using this and the open-commutativity
of the functor $\lambda$ we can in the same way conclude that the
functor $\lambda_{cc}$ is also open-multicommutative.

Q.E.D.

\medskip

Let us consider now the functors $G_{cc}P$ and $\lambda_{cc}P$ which are defined
in Teleiko and Zarichnyi (1999). Actually,
$G_{cc}P\left(X\right)=G_{cc}\left(P\left(X\right)\right)$ for every compactum
$X$ and it is the composition of two functors $G_{cc}$ and $P$ (the same
situation is for the functor $\lambda_{cc}P$). The following result can be
derived from Propositions 4 and 6.

{\bf Corollary 2.} {\it The functors $G_{cc}P$ and $\lambda_{cc}P$
are open-multicommutative.}

\end{document}